\documentclass[10pt,a4paper,twoside]{amsart} 

\usepackage{amsfonts}
\usepackage{amsmath}
\usepackage{amsthm}
\usepackage{amscd}
\usepackage{euscript}

\DeclareMathOperator{\End}{End}
\DeclareMathOperator{\Ext}{Ext}
\DeclareMathOperator{\HH}{H}
\DeclareMathOperator{\Hom}{Hom}
\DeclareMathOperator{\codim}{codim}
\DeclareMathOperator{\coker}{coker}

\DeclareMathOperator{\reg}{reg}

\newcommand{\EE}{{\EuScript E}}
\newcommand{\FF}{{\EuScript F}}
\newcommand{\II}{{\mathfrak I}}
\newcommand{\MM}{{\mathfrak M}}
\newcommand{\OO}{{\EuScript O}}
\newcommand{\PPii}{{{\mathbb P}^2}}
\newcommand{\PPn}{{{\mathbb P}^n}}
\newcommand{\PP}{{\mathbb P}}

\newcommand{\ZZ}{{\mathbb Z}}

\newcommand{\dual}[1]{{#1}^*}
\newcommand{\second}{{\prime\prime}}

\swapnumbers
\newtheorem{thm}{Theorem}[section]
\newtheorem{lemma}[thm]{Lemma}
\newtheorem{cor}[thm]{Corollary}
\newtheorem{prop}[thm]{Proposition}
\theoremstyle{definition}
\newtheorem{defin}[thm]{Definition}

\theoremstyle{remark}

\newtheorem{exm}[thm]{Example}

\begin{document}
\title{Minimal resolution of general stable vector bundles on $\PP^2$}

\author{Carla Dionisi and Marco Maggesi}

\address{Carla Dionisi \\
  Dipartimento di matematica e applicazioni R.~Caccioppoli \\
  via Cintia loc.\ Monte S.~Angelo \\
  80138 Napoli, Italy}

\email{dionisi@matna2.dma.unina.it}

\address{Marco Maggesi \\
  Dipartimento di Matematica U.~Dini \\
  Viale Morgagni 67/a \\
  50134 Firenze, Italy}

\email{maggesi@math.unifi.it}


\begin{abstract}
  We study the general elements of the moduli spaces \( \MM_{\PP^2}
  (r, c_1, c_2) \) of stable holomorphic vector bundle on $\PP^2$ and
  their minimal free resolution.  Incidentally, a quite easy proof of
  the irreducibility of \( \MM_{\PP^2} (r, c_1, c_2) \) is shown.
\end{abstract}

\maketitle

\section{Introduction} \label{sec:intro}
In this paper we investigate stable vector bundles on the complex
projective plane $\PP^2$ by means of their minimal free resolution.
The fundamental background for this study is the work of Bohnhorst and
Spindler \cite{bs} and their idea to use admissible pairs to
characterize the stability of rank-$n$ vector bundles on $\PP^n$ and
to give a stratification of the relative moduli space in constructible
subsets.

Two main difficulties arise: (i) to state a weak version of the
Bohnhorst-Spindler theorems for rank $r \geq 2$ vector bundles on
$\PP^2$, (ii) to estimate the codimension of the constructible subsets
of the moduli space.  In section $2$ we make some general remarks on
rank $r$ vector bundles on $\PP^n$ with $r \geq n$ to address (i),
while (ii) is the object of the last section.

We would like to thank G.~Ottaviani for his invaluable guidance and
V.~Ancona for many useful discussions.

\section{Admissible pairs and resolutions}

Let 
\begin{equation}
  \label{res:bs}
  0 \to \bigoplus_{i=1}^k \OO_\PPn (-a_i) \xrightarrow{\Phi}
  \bigoplus_{j=1}^{r+k} \OO_\PPn (-b_j) \to \EE \to 0
\end{equation}
be a free resolution of length $1$ of a rank $r$ vector bundle $\EE$
on $\PP^n$.  We assume that the two sequences $a_i$ and $b_i$ are
indexed in nondecreasing order
\begin{equation}
  \label{eq:aibj}
  \begin{aligned}
    a_1 &\leq a_2 \leq \dots \leq a_k, \\
    b_1 &\leq b_2 \leq \dots \leq b_k \leq \dots \leq b_{r+k}.
  \end{aligned}
\end{equation}
We call \( (a,b) = ((a_1, \dots, a_k), (b_1, \dots, b_{r+k})) \) the
\emph{pair} associated to the resolution (\ref{res:bs}).  If the
resolution (\ref{res:bs}) is minimal, then we call $(a,b)$ the pair
associate to the bundle $\EE$.  Notice that the associated pair and
the Betti numbers of a resolution encode exactly the same information;
in particular \( \max(a_k - 1, b_{r+k}) \) is the regularity.

\begin{defin}
  The pair $(a,b)$ is said to be weakly admissible if
  \begin{equation}
    \label{eq:wadmissible}
    a_i > b_{n+i} \qquad\text{for all $i=1,\dots,k$}
  \end{equation}
  and admissible (or strongly admissible) if
  \begin{equation}
    a_i > b_{r+i} \qquad\text{for all $i=1,\dots,k$}
  \end{equation}
  For brevity we say that the resolution (\ref{res:bs}) is weakly or
  strongly admissible if the associated pair $(a,b)$ is.
\end{defin}

\begin{exm}
  \label{exm:banded}
  The map $\Phi$ can be expressed by a \( (r+k)\times k \) matrix of
  forms \( (\phi_{i, j}) \) of degree \( \deg (\phi_{i, j}) = (b_i -
  a_j) \).  If $(a,b)$ is a strongly admissible pair and \( \omega_0
  \ldots \omega_r \) are linear forms in general position on $\PP^n$,
  then the \( (r+k) \times k \) matrix
  \begin{equation}
    (\phi_{ij}) :=
    \begin{bmatrix}
      \omega_0^{a_1-b_1}     & \dots  & 0                        \\
      \vdots                 & \ddots & \vdots                   \\
      \vdots                 &        & \omega_0^{a_k-b_k}       \\
      \omega_r^{a_1-b_{r+1}} &        & \vdots                   \\
      \vdots                 & \ddots & \vdots                   \\
      0                      & \dots  & \omega_r^{a_k-b_{r+k}}
    \end{bmatrix}
  \end{equation}
  defines a minimal free resolution and \( \EE := \coker \Phi \) is a
  vector bundle with associated pair $(a,b)$.  If the pair $(a,b)$ is
  only weakly admissible, the same reasoning works for the $\omega$'s
  defined by
  \begin{equation*}
    (\omega_0,\, \dots,\, \omega_n,\,\omega_{n+1},\, \dots,\, \omega_r)
    = (x_0,\, \dots,\, x_n,\, 0,\, \dots,\, 0)
  \end{equation*}
\end{exm}

Admissibility was originally introduced by Bohnhorst and Spindler in
\cite{bs} to characterize stability for vector bundles on $\PP^n$ of
homological dimension $1$ with rank $r$ equal to $n$.  Note that for
$r=n$ weakly and strongly admissibility coincide.  In this section we
are going to restate their results with some generalizations to the
case of rank $r\geq n$.

\begin{prop}
  If $r\geq n$, the following two condition on the resolution
  (\ref{res:bs}) are equivalent:
  \begin{enumerate}
  \item is minimal;
  \item is weakly admissible and every constant entry of the matrix
    $(\phi_{i,j})$ is zero.
  \end{enumerate}
\end{prop}
\begin{proof}
  Obviously, (2.)\ follows from (1.).  For $r=n$ the statement was
  proved by Bohnhorst and Spindler (\cite{bs} proposition 2.3).  Now
  suppose that $r>n$ and (\ref{res:bs}) is minimal.  Since
  $\EE(b_{r+k})$ is globally generated, Bertini's theorem ensures that
  a generic map \(f \colon \OO(-b_{r+k}) \to \EE\) is injective.
  Then, in the following commutative diagram columns and rows are
  exact and $\EE^{\second}$ is locally free:
  \begin{small}
    \begin{equation}
      \begin{CD}
        {} @. 0 @. 0 @. 0 @. {} \\
        @. @VVV @VVV @VVV @. \\
        0 @>>> 0 @>>> \OO(-b_{r+k}) @>Id>> \OO(-b_{r+k}) @>>> 0 \\
        @. @VVV @VVV @VVV @. \\
        0 @>>> 
        \oplus_{i=1}^k \OO(-a_i) @>>>
        \oplus_{i=1}^{r+k} \OO(-b_i) @>>>
        \EE @>>> 0 \\
        @. @VV{Id}V @VVV @VVV @. \\
        0 @>>> 
        \oplus_{i=1}^k \OO(-a_i) @>>>
        \oplus_{i=1}^{r+k-1} \OO(-b_i) @>>>
        \EE^{\second} @>>> 0 \\
        @. @VVV @VVV @VVV @. \\
        {} @. 0 @. 0 @. 0 @. {}
      \end{CD}
    \end{equation}
  \end{small}
  The minimality of the middle row yields the minimality of the last
  row.  Using induction on $r$, we may assume that the last row is
  weakly  admissible.  Then, the middle row is also weakly admissible.
\end{proof}

As a consequence, to every vector bundle it corresponds a weakly
admissible pair.  Vice versa, for every weakly admissible pair, the
examples \ref{exm:banded} provide a way to construct a vector bundle
with such associated pair.

\begin{thm}[Bohnhorst-Spindler \cite{bs}]
  \label{th:bs}
  Suppose $r=n$ and that the resolution (\ref{res:bs}) is admissible.
  Let $c_1=\sum a_i - \sum b_j$ be the first Chern class and
  $\mu=c_1/n$ the slope of $\EE$.  Then $\EE$ is stable (respectively
  semistable) if and only if
  \begin{equation}
    \label{eq:semistable}
    b_1 > - \mu \qquad \text{(resp. $b_1 \geq - \mu$)}.
  \end{equation}
\end{thm}

For the case of higher rank, we have no chances to extend the above
arithmetical characterization since stable and unstable vector bundles
may have the same associated pair.  However, one implication still
hold:

\begin{thm}
  \label{thm:rgtn}
  If the resolution (\ref{res:bs}) is weakly admissible (in particular
  if it is minimal) and $\EE$ is a stable (resp. semistable) vector
  bundle, then the associated pair $(a,b)$ is strongly admissible and
  \(b_1> -\mu\) (resp. \(b_1 \geq -\mu\)).
\end{thm}

\begin{proof}
  We first prove that if $\EE$ is semistable, then $b_{r+k}<a_k$.
  In fact, if $b_{r+k}\geq a_k$ then $\EE$ split as \(\EE =
  \OO(-b_{k+r}) \oplus \EE^{\second}\) and by weakly admissibility
  \begin{equation}
    \begin{split}
      \sum_{i=1}^k a_i - \!\!\sum_{j=1}^{r+k-1}b_j
      &= - \sum_{i=1}^n b_i + \sum_{i=1}^k (a_i - b_{n+i}) -
      \!\!\!\!\sum_{i=n+k+1}^{r+k-1} b_i \\
      &\geq -n b_{r+k} + k - (r - n - 1) b_{r+k} \\
      &> (1-r)b_{r+k}
    \end{split}
  \end{equation}
  then we have
  \begin{math}
    \mu(\EE^{\second}) > - b_{r + k} = \mu(\OO(- b_{r + k}))
  \end{math}
  which contradicts the semistability of $\EE$.
  
  Now suppose that $\EE$ is semistable and $a_s\leq b_{r+s}$ for some
  $s$ with $1 \leq s < k$ and let $(s_0-1)$ be the largest of such
  $s$.  Since 
  \begin{equation*}
    a_1\geq \dots \geq a_{s_0} \geq b_{r+s_0} \geq \dots \geq b_{r+k},
  \end{equation*}
  the minor $\Phi^{\second}$ of $\Phi$, obtained by cutting off the
  last $(k-s_0)$ rows and $(k-s_0)$ columns, remains of maximal rank
  so \( \EE^\second := \coker \Phi^\second \) is a vector bundle.  A
  surjective morphism $\EE \to \EE^{\second} \to 0$ is defined by the
  diagram
  \begin{small}
    \begin{equation}
      \begin{CD}
        \label{eq:destabilizza}
        0 @>>>
        \oplus_{i=1}^{k} \OO(-a_i) @>{\Phi}>>
        \oplus_{j=1}^{r+k} \OO(-b_j) @>>>
        \EE @>>> 0 \\
        @. @VVV @VVV @VVV @. \\
        0 @>>>
        \oplus_{i=1}^{s_0} \OO(-a_i) @>{\Phi^{\second}}>>
        \oplus_{j=1}^{r+s_0}\OO(-b_j) @>>>
        \EE^{\second} @>>> 0
      \end{CD}
    \end{equation}
  \end{small}
  where the first two vertical map are the natural projections.
  Observe that $\EE^{\second}$ have the same rank $r$ as $\EE$ and
  \begin{equation}
    \begin{split}
      \mu(\EE^{\second}) &=
      \frac{1}{r} \left(
        \sum_{i=1}^{s_0} a_i - \sum_{j=1}^{r+s_0} b_j
      \right)
      = \mu(\EE) - \frac{1}{r} \sum_{i=s_0+1}^k (a_i - b_{r+i}) \\
      &< \mu(\EE).
    \end{split}
  \end{equation}
  Then, $\EE^{\second}$ must be semistable otherwise any torsionless
  quotient sheaf destabilizing it would also destabilize $\EE$.  By
  induction on $k$, we may assume that the second row of
  (\ref{eq:destabilizza}) is strongly admissible.  In particular, we
  have \( a_{s_0} > b_{r+s_0} \) that gives a contradiction.
  
  Finally, if $\EE$ is stable (resp. semistable), then
  \begin{equation}
    \HH^0(\EE(m)) = 0 \qquad \forall m \leq -\mu(\EE)
    \quad \text{(resp. \(\forall m < -\mu(\EE)\))}
  \end{equation}
  but, from the exact sequence
  \begin{equation}
    0 \to
    \bigoplus_{i=1}^k \OO(-a_i+b_1) \to
    \bigoplus_{j=1}^{r+k} \OO(-b_j+b_1) \to
    \EE(b_1) \to 0,
  \end{equation}
  we have \(\HH^0(\EE(b_1)) \neq 0\) then $b_1 > -\mu(\EE)$ (resp.
  $b_1 \geq -\mu(\EE)$).
\end{proof}

\section{On vector bundles on $\PP^2$} \label{sec:peetwo}
By Horrocks theorem, every rank $r$ vector bundle $\EE$ on $\PP^2$ has
homological dimension at most $1$, that is, if $\EE$ does not split in
the direct sum of line bundles, then it is presented by a minimal free
resolution of the form
\begin{equation}
  \label{res:bs2}
  0 \to \bigoplus_{i=1}^k \OO_\PPii (-a_i) \xrightarrow{\Phi}
  \bigoplus_{j=1}^{r+k} \OO_\PPii (-b_j) \to \EE \to 0.
\end{equation}

The Chern classes $c_1$, $c_2$ of $\EE$ are determined by
$a_i$ and $b_j$ with the formulas
\begin{equation}
  \begin{aligned}
    \label{eq:thechernclasses}
    c_1 &= \sum_{i=1}^k a_i - \sum_{i=1}^{k+r} b_i, \\
    2c_2 - c_1^2 &= \sum_{i=1}^k a_i^2 - \sum_{i=1}^{k+r} b_i^2.
  \end{aligned}
\end{equation}

We denote by $\II$ the set of all (strongly) admissible pairs $(a,b)$
associated to rank $r$-vector bundles on $\PP^2$ with Chern classes
$c_1$, $c_2$ satisfying the condition \( b_1 > -\mu = (\sum a_i - \sum
b_j) / r \).  Theorem \ref{thm:rgtn} shows that the set $\II$ contains
the set of all possible associated pair to a stable vector bundle in
\( \MM_\PPii (r, c_1, c_2) \) and coincides exactly with it for $r=n$.
Then
\begin{equation}
  \MM_\PPii (r, c_1, c_2) = \coprod_{(a, b) \in \II} \MM (a, b)
\end{equation}
where \( \MM (a, b) \) will be the subset (possibly empty) of \(
\MM_\PPii (r, c_1, c_2) \) of vector bundles with associated pair
$(a,b)$.

The following result was stated and proved by Bohnhorst and Spindler
\cite{bs} for rank-$n$ vector bundles on $\PP^n$ with homological
dimension $1$, but their proof works on $\PP^2$ for vector bundles of
any rank without modifications.

\begin{thm}
  \label{th:codim}
  For all \( (a, b) \in \II \), the closed set \( \overline{\MM(a, b)}
  \) is an irreducible algebraic subset of \( \MM_\PPii (r, c_1, c_2)
  \) of dimension:
  \begin{equation}
    \label{eq:dim}
    \begin{split}
      \dim\overline{\MM(a,b)}
      &= \dim\Hom(F_1,F_0) + \dim\Hom(F_0,F_1) \\
      & \quad - \dim\End(F_1) - \dim\End(F_0) 
      + 1 - \# \{ (i,j) : a_i=b_j \},
    \end{split}
  \end{equation}
  where
  \begin{math}
    F_0 = \oplus_{j=1}^{k+r}\OO(- b_j),
  \end{math}
  \begin{math}
    F_1 = \oplus_{i=1}^k \OO(- a_i).
  \end{math}
\end{thm}

\section{Natural pairs and general vector bundles}
We say that 
\begin{math}
  (a, b) = ((a_1, \dots, a_k), (b_1, \dots, b_{r+k}))
\end{math}
is a \emph{natural pair} if it is admissible and
\begin{equation}
  \label{eq:natural}
  b_{r+k} < a_1, \qquad a_k \leq b_1 + 2.
\end{equation}
Through this section, we are going to show that resolutions of general
vector bundles have natural pairs:

\begin{thm}
  \label{thm:null-codim}
  One has \( \codim \overline{\MM(a, b)} = 0 \) if and only if $(a,b)$
  is a natural pair.
\end{thm}

As a remarkable consequence we will derive a quite simple proof of the
irreducibility of the moduli spaces of stable vector bundles on
$\PP^2$ (other proofs with different techniques can be found in
\cite{ba1}, \cite{ell}), \cite{hl}, \cite{LeP}, \cite{M}) and we will
compute the regularity and the cohomology of their general elements.

We recall that, since \( \dim \Ext^2(\FF, \FF) = 0 \) for any stable
vector bundle $\FF$ on $\PP^2$, the relative moduli \( \MM_\PPii(r,
c_1, c_2) \) space is smooth of dimension
\begin{equation}
  \label{eq:dimension}
  \dim \Ext^1(\FF, \FF) = 2 r c_2 - (r - 1) c_1^2 - r^2 + 1.  
\end{equation}

Let us consider the function \( A(t) := h^2(\OO(t)) \) and the finite
differences of first and second order
\begin{math}
  (\Delta_u A)(t) := A(t + u) - A(t)
\end{math}
and
\begin{math}
  (\Delta_v \Delta_u A)(t) := (\Delta_u A)(t + v) - (\Delta_u A)(t).
\end{math}

\begin{lemma}
  \label{lemma:codimension}
  Let $(a,b)$ be the admissible pair associated to a stable vector
  bundle $\EE$ on $\PP^2$.  Then
  \begin{equation}
    \label{eq:codimension}
    \begin{split}
      \codim \overline{\MM(a,b)}
      &= \sum_{i=1}^r h^1( \EE(b_i)) + \# \{ (i,j) : a_i=b_j \} \\
      &\quad + \sum_{i,j=1}^r (\Delta_{b_{i+r}-a_i}
      \Delta_{b_{j+r}-a_j} A) (a_i - b_{j+r})
    \end{split}
  \end{equation}
\end{lemma}

\begin{proof}
  Let 
  \begin{equation}
    \label{eq:f0f1}
    0 \to F_1 \to F_0 \to \EE \to 0
  \end{equation}
  be the minimal resolution of $\EE$ where
  \begin{equation}
    F_0 = \bigoplus_{j=1}^{k+r}\OO(- b_j), \qquad
    F_1 = \bigoplus_{i=1}^k \OO(- a_i).
  \end{equation}
  
  The stability of $\EE$ ensures the vanishing \( \dim (\Ext^2 (\EE,
  \EE)) = h^2 (\dual {\EE} \otimes \EE) = 0 \) so that \(
  h^2(\dual{F_0} \otimes \EE ) = h^2(\dual{F_1} \otimes \EE) \).
  Then, from (\ref{eq:f0f1}) we easily find the following data:
  \begin{equation}
    \label{eq:dati}
    \begin{split}
      h^0(\dual{F_0} \otimes \EE)
      &=h^0(\dual{F_0} \otimes F_0)-h^0(\dual{F_0} \otimes F_1), \\
      h^0(\dual{F_1} \otimes \EE)
      &=h^0(\dual{F_1} \otimes F_0)-h^0(\dual{F_1} \otimes F_1), \\   
      \dim(\Ext^1(\EE,\EE)) &= h^1(\dual{\EE} \otimes \EE) = \\
      &=h^1(\dual{F_0} \otimes \EE )-h^1(\dual{F_1} \otimes \EE )+ \\
      &\quad + h^0(\dual{F_1} \otimes \EE )- h^0(\dual{F_0} \otimes \EE)+1
    \end{split}
  \end{equation}            
  and from (\ref{eq:dim}) we have
  \begin{equation}
    \label{eq:codim1}
    \begin{split}
      \codim \overline{\MM(a,b)}
      &= \dim(\Ext^1(\EE,\EE)) - \dim \overline{\MM(a,b)}= \\
      &= h^1(\dual{F_0} \otimes \EE ) - h^1(\dual{F_1} \otimes \EE)
      + \# \{ (i,j) : a_i=b_j \}
    \end{split}
  \end{equation}
  Now, by splitting $F_0$ as
  \begin{math}
    \OO(-b_1) \oplus \dots \oplus \OO(-b_{r-1})
    \oplus \tilde{F_0}
  \end{math}
  with \( \tilde {F_0} := \oplus_{i = r + 1}^{k + r}\OO(- b_i) \), the
  above formula becomes
  \begin{equation}
    \label{eq:codim}
    \begin{split}
      \codim \overline{\MM(a,b)}
      &= \sum_{i=1}^r h^1( \EE(b_i))
      + \# \{ (i,j) : a_i=b_j \} \\
      &\quad + h^1(\dual{\tilde{F_0}} \otimes \EE )
      - h^1(\dual{F_1} \otimes \EE) \\
      &= \sum_{i = 1}^r h^1(\EE(b_i))
      + \# \{ (i,j) : a_i=b_j \} \\
      &\quad + h^2(\dual{\tilde{F_0}} \otimes F_1)
      - h^2(\dual{\tilde{F_0}} \otimes F_0) \\
      &\quad - h^2(\dual{F_1} \otimes F_1 )
      + h^2(\dual{F_1} \otimes F_0).
    \end{split}
  \end{equation}
  Since
  \begin{math}
    h^2(\dual {\tilde {F_0}} \otimes F_0) =
    h^2(\dual{\tilde{F_0}} \otimes \tilde{F_0} )
  \end{math}
  and
  \begin{math}
    h^2(\dual{F_1} \otimes F_0) =
    h^2(\dual{F_1} \otimes \tilde{F_0})
  \end{math}
  so
  \begin{equation}
    \label{eq:CODIME}
    \begin{split}
      \codim \overline{\MM(a,b)}
      &= \sum_{i=1}^r h^1( \EE(b_i))
      + \# \{ (i,j) : a_i=b_j \} \\
      &\quad + h^2(\dual{\tilde{F_0}} \otimes F_1)
      -h^2(\dual{\tilde{F_0}} \otimes \tilde{F_0}) \\
      &\quad - h^2(\dual{F_1} \otimes F_1 )
      + h^2(\dual{F_1} \otimes \tilde{F_0}).
    \end{split}
  \end{equation}
  Finally, distributing the direct sums appearing in the definition of
  $\tilde F_0$ and $F_1$ the equation (\ref{eq:CODIME}) becomes
  \begin{equation}
    \begin{split}
      \codim \overline{\MM(a,b)}
      &= \sum_{i=1}^r h^1( \EE(b_i))
      + \# \{ (i,j) : a_i=b_j \} \\
      &\quad + \sum_{i,j=1}^r
      \Big[
      h^2 (\OO (b_{i+r} - a_j)) - h^2 (\OO (b_{i+r} - b_{j+r})) \\
      &\qquad\qquad - h^2 (\OO (a_i - a_j)) + h^2 (\OO (a_i - b_{j+r}))
      \Big] \\
      &= \sum_{i=1}^r h^1( \EE(b_i))
      + \# \{ (i,j) : a_i=b_j \} \\
      &\quad + \sum_{i,j=1}^r
      (\Delta_{b_{i+r}-a_i} \Delta_{b_{j+r}-a_j} A)
      (a_i - b_{j+r})
    \end{split}
  \end{equation}
\end{proof}

We observe that natural pairs are parametrized by
three integers $s$, $k$, $\alpha$ such that
\begin{equation}
  \label{eq:skalpha}
  k \geq 1 \qquad \text{and} \qquad-k+1\leq \alpha \leq k+r
\end{equation}
as follows: the pair \( (a, b)_{s, k, \alpha} \) corresponding to the
triple \( (s, k, \alpha) \) is the pair associated to a resolution of
the form
\begin{equation}
  \label{eq:res:alphapositive}
  0 \to
  \OO(-s-1)^k \to
  \OO(-s)^\alpha \oplus \OO(-s+1)^{r+k-\alpha} \to
  \EE \to 0
\end{equation}
if $\alpha \geq 0$, or of the form
\begin{equation}
\label{eq:res:alphanegative}
  0 \to \OO(-s-1)^{k+\alpha} \oplus \OO(-s)^{-\alpha} \to
  \OO(-s+1)^{r+k} \to
  \EE \to 0
\end{equation}
if $\alpha < 0$.  We exclude the case $\alpha=-k$ so that $s$ is the
regularity of the pair, i.e.\ \( s = \max (a_k - 1, b_{r+k}) \).

\begin{proof}[Proof of theorem \ref{thm:null-codim}]
  It can be verified by direct computation from theorem \ref{th:codim}
  that, if $\EE$ has natural pair, then the codimension of \(
  \overline {\MM(a, b)} \) is zero.  Conversely, let $u$, $v$ be two
  non-negative integers.  Since all finite difference \( (\Delta_u
  A)(t) := A(t + u) - A(t) \) are non decreasing functions of $t$,
  then
  \begin{equation}
    \label{eq:delta}
    (\Delta_v \Delta_u A) (t) \geq 0 
  \end{equation}
  and by the previous lemma
  \begin{equation}
    \label{eq:positivo}
    \codim \overline{\MM(a,b)} \geq
    \sum_{i=1}^r h^1( \EE(b_i)) +
    \# \{ (i,j) : a_i=b_j \}.
  \end{equation}
  If \( \codim \overline{\MM(a, b)} = 0 \), we have \( a_k \leq b_1 +
  2 \) and \( \# \{ (i,j) : a_i=b_j \} = 0 \), since \( h^1 (\EE
  (b_1)) = 0 \) implies \( h^2 (F_1 (b_1)) = 0 \).  This forces
  $(a,b)$ to be a natural pair.
\end{proof}

\begin{prop}
  \label{prop:irr}
  Let \( \MM_{\PP^2} (r, c_1, c_2) \) be nonempty and
  \begin{equation}
    \label{eq:esse}
    s := \max \{ \rho \in \ZZ :
    r{\rho}^2 + 2 c_1 \rho - r \rho \leq 2 c_2 - c_1^2 + c_1 - 1 \},
  \end{equation}
  or, equivalently,
  \begin{equation}
    \label{eq:esse2}
    s := \min\{\rho\in\ZZ : r\rho^2 + 2c_1\rho + r\rho \geq
    2c_2 -c_1^2 -c_1 \}.
    \tag{\ref{eq:esse}bis}
  \end{equation}
  If $\alpha$ and $k$ are defined by
  \begin{equation}
    \label{eq:alphak}
    \begin{aligned}
      \alpha &:= 2 c_2 - c_1^2 + r - r s^2 - 2 c_1 s, \\ 
      k & := (r s + c_1 - r + |\alpha|) / 2,
    \end{aligned}
  \end{equation}
  then \( (a, b)_{s,k,\alpha} \) is the only natural pair of \(
  \MM_{\PP^2} (r, c_1, c_2) \).
\end{prop}

\begin{proof}
  This is a verification; we outline the main steps of the
  computation.  In the first place, one must ensure that the natural
  pair \( (a,b)_{s,k,\alpha} \) is actually associated to vector
  bundles in \( \MM_{\PP^2} (r, c_1, c_2) \).  This amount to show
  that from (\ref{eq:thechernclasses}) the pair \( (a,b)_{s,k,\alpha}
  \) has the appropriate Chern classes and that conditions
  (\ref{eq:skalpha}) hold.
  
  From theorem \ref{thm:null-codim}, any pair $(a,b)$ such that \(
  \dim \overline {\MM (a, b)} = 0 \) is a natural pair of the form \(
  (a, b)_{s,k,\alpha} \).  From resolutions
  (\ref{eq:res:alphapositive}) and (\ref{eq:res:alphanegative}) we
  find that $\alpha$, $k$ must satisfy (\ref{eq:alphak}).  Then, it
  remains to verify that $s$ is uniquely determined from $r$, $c_1$,
  $c_2$ and satisfy (\ref{eq:esse}).  By substitution, the
  inequalities \( -k < \alpha \leq k + r \) turn into
  \begin{equation}
    r s^2 + 2 c_1 s - c_1 - r s + 1 \leq
    2 c_2 - c_1^2 \leq
    r s^2 + 2 c_1 s + c_1 + r s.
  \end{equation}
  Since the intervals
  \begin{math}
    [ r s^2 + 2 c_1 s - c_1 - r s + 1,
    \quad r s^2 + 2 c_1 s + c_1 + r s ]
  \end{math}
  are disjoint for $s$ varying in $\ZZ$, then equations
  (\ref{eq:esse}) and (\ref{eq:esse2}) give the only suitable value
  for $s$.
\end{proof}

\begin{thm}
  \label{thm:irr}
  The moduli spaces of stable rank $r$ vector bundles on $\PP^2$ are
  irreducible.
\end{thm}

\begin{proof}
  Moduli space of stable rank $r$ vector bundles on $\PP^2$ are
  smooth.  By the previous proposition they can have only one
  connected component.
\end{proof}

\begin{cor}
  The general element of \( \MM_{\PP^2}(r,c_1,c_2) \) has natural
  cohomology.
\end{cor}

The above corollary justify the terminology \emph{``natural pair''}.
A different proof for it can be found in \cite{hl}.

Now, we are going to give some inequalities on the regularity and the
cohomology of stable vector bundles using proposition \ref{prop:irr}.
In particular, for rank $2$ vector bundles, the next two corollaries
give respectively a refined version of corollary 5.4 \cite{br} and
proposition 7.1 \cite{ha78}.

\begin{cor}
  A general vector bundle $\EE$ in \( \MM_{\PP^2} (r, c_1, c_2) \) has
  regularity \( \reg (\EE) = s \), where $s$ is given by
  (\ref{eq:esse}).
\end{cor}

\begin{cor}
  \label{cor:Ha}
  Let $[\EE]$ be a vector bundle in \( \MM = \MM_{\PP^2} (r, c_1, c_2)
  \) and $s$ defined by (\ref{eq:esse}). Then \( \HH^0(\EE(t)) \neq
  0\) if
  \begin{equation*}
    \begin{aligned}
      t&\geq s   &&\quad \text{when} \quad rs^2+2c_1s+rs = 2c_2-c_1^2-c_1, \\
      t&\geq s-1 &&\quad \text{otherwise.}
    \end{aligned}
  \end{equation*}
  The above inequality is sharp, in the sense that, if $\EE$ is
  general, it gives a necessary and sufficient condition.
\end{cor}

\begin{proof}
  Let \( ((a_1,\dots,a_k), (b_1,\dots,b_{k+r})) \) be the admissible
  pair associated to a vector bundle $\EE$ in $\MM$.  Then, one has
  \(\HH^0(\EE(t)) \neq 0\) if and only if $t - b_1 \geq 0$.
  By semicontinuity of cohomology groups and theorem \ref{thm:irr}, it
  is enough to restrict ourselves to the case where $\EE$ is general.
  So, by (\ref{eq:res:alphapositive}) and (\ref{eq:res:alphanegative})
  one has \(\HH^0(\EE(t)) \neq 0\) if and only if 
  \begin{equation*}
    \begin{aligned}
      t&\geq s && \quad \text{if} \quad \alpha=k+r\\
      t&\geq s-1 && \quad \text{otherwise}
    \end{aligned}
  \end{equation*}
  and the condition $\alpha=k+r$ is equivalent to \( rs^2 + 2c_1s + rs
  = 2c_2 - c_1^2 - c_1\) by (\ref{eq:alphak}).
\end{proof}


\providecommand{\bysame}{\leavevmode\hbox to3em{\hrulefill}\thinspace}

\end{document}